\documentclass[11pt]{article}

\usepackage{latexsym}
\usepackage{amssymb,hyperref}
\newtheorem{theoreme}{{\bf Th\'eor\`eme}}[section]

\newtheorem{corollaire principal}[principal]{\bf Corollaire}
\newtheorem{proposition}[theoreme]{{\bf Proposition}}

\newtheorem{corollaire}[theoreme]{{\bf Corollaire}}

\newenvironment{demonstration}{\noindent{\bf D\'emonstration
}}{\nolinebreak $\Box $\hspace{-2.15mm}\raisebox{1.25mm}{.} \medskip}

\newenvironment{demonstration du lemme}{\noindent{\bf D\'emonstration du lemme
}}{\nolinebreak $\Box $\hspace{-2.15mm}\raisebox{1.25mm}{.} \medskip}

\def\CC{{\bf C}}

\def\QQ{{\bf Q}}

\newtheorem{theo}{Theorem}[section]
\newtheorem{prop}[theo]{Proposition}
\newtheorem{cor}[theo]{Corollary}

\begin{document}

\title{{\Large \bf  Une caract\'erisation des vari\'et\'es complexes  compactes parall\'elisables  admettant des structures affines}}

\author{{\normalsize
{\bf Sorin DUMITRESCU}}}
\date{Janvier 2009}

\maketitle

\vspace{0.1cm}

\noindent{\bf{ R\'esum\'e.}} 
Nous caract\'erisons les vari\'et\'es complexes compactes parall\'elisables  admettant des connexions affines holomorphes plates sans torsion. En particulier, nous exhibons des vari\'et\'es complexes compactes admettant des connexions affines holomorphes, mais aucune connexion affine holomorphe plate sans torsion.

\noindent{\bf{ Abstract.}}  {\bf On complex  compact parallelizable manifolds admitting   affine structures.}
We classify complex compact parallelizable manifolds which admit flat torsion free holomorphic affine connections. We exhibit complex compact manifolds admitting holomorphic affine
connections, but no flat torsion free holomorphic affine connections.

{\bf Abridged English Version}

We deal with compact complex manifolds admitting holomorphic affine connections.
The main results  are  the following :

\begin{prop} Let $M$ be a complex compact  parallelizable manifold which is a quotient  $G/ \Gamma$ of a $n$-dimensional complex connected  Lie group  $G$ by an uniform lattice $\Gamma$. Then:

(i) $M$ admits holomorphic affine connections. The pull-back of such a connection on $G$ is invariant by right translation.

(ii) $M$ admits flat holomorphic affine connections.

(iii) $M$ admits flat torsion free holomorphic affine connections if and only if there is a Lie algebra monomorphism $i : \mathcal G \to gl(n, \CC) \ltimes \CC^n$ such that
$i(\mathcal G)$ intersects trivially the isotropy $gl(n,\CC)$, where $\mathcal G$ is the Lie algebra of $G$ and $gl(n, \CC) \ltimes \CC^n$ is the Lie algebra of the complex affine groupe
of $\CC^n$.
\end{prop}

\begin{cor} If $G$ is a complex semi-simple Lie group, then a compact  complex parallelizable  manifold $M= G / \Gamma$ doesn't admit flat torsion free holomorphic affine connections.
\end{cor}

\begin{cor} Let $M=G / \Gamma$ be a three dimensional compact complex parallelizable manifold . Then $M$ admits flat torsion free holomorphic affine connections
if and only if $G$ is solvable.
\end{cor}

 We also exhibit examples of non parallelizable compact complex
manifolds which admit (flat) holomorphic affine connections, but no flat torsion free holomorphic affine connections.

The natural question which is still open is whether  a compact complex manifold admitting  holomorphic affine connections, also  admit a flat holomorphic affine connection.

{\it  Acknowledgement. This work was  partially supported by  the ANR Grant Symplexe BLAN 06-3-137237.}

\section{Introduction}

Nous nous int\'eressons ici aux vari\'et\'es complexes compactes $M$ dont le fibr\'e tangent holomorphe $TM$ admet des connexions affines holomorphes. 

Rappelons qu'une connexion affine holomorphe $\nabla$
est plate (de courbure nulle)  et sans torsion si et seulement si elle est  localement isomorphe \`a la connexion standard de $\CC^n$. Dans ce cas  $\nabla$ munit la vari\'et\'e $M$
d'une {\it structure affine complexe}, autrement dit $M$ admet un atlas compatible avec la structure complexe dont les applications de changement de carte sont des transfomations
affines (voir, par exemple,~\cite{Gunning,Kobayashi}). R\'eciproquement, une structure affine complexe sur $M$ induit canoniquement une connexion affine holomorphe plate et
sans torsion.

Si la vari\'et\'e $M$ est
suppos\'ee k\"ahl\'erienne, la th\'eorie de Chern-Weil assure que les classes de Chern \`a coefficients rationnels $c_{i}(M, \QQ)$ s'annulent et il n'existe aucune obstruction
topologique \`a l'existence d'une connexion affine holomorphe  plate sur $M$~\cite{Gunning}.  Effectivement,  il est montr\'e dans $\cite{IKO}$  que $M$ admet un rev\^etement fini non ramifi\'e qui est un tore complexe, quotient de $\CC^n$ par un r\'eseau de translations. En particulier, $M$ h\'erite
de la structure affine complexe  de $\CC^n$ (donc d'une connexion affine holomorphe plate sans torsion). 

Il est \'egalement prouv\'e dans $\cite{IKO}$ que {\it toute surface complexe compacte admettant une connexion affine holomorphe  est biholomorphe   \`a un quotient d'un ouvert de $\CC^2$ par un sous-groupe discret de transformations affines agissant proprement et sans point fixe.} En particulier, $M$ admet une connexion 
affine holomorphe plate sans torsion.

Le but de cette note est de construire des vari\'et\'es complexes compactes admettant des connexions affines holomorphes, mais aucune connexion affine holomorphe plate sans torsion. Pour cela nous caract\'erisons les vari\'et\'es holomorphes parall\'elisables qui admettent des structures affines complexes. Rappelons que, d'apr\`es~\cite{Wa}, une vari\'et\'e  complexe compacte  parall\'elisable
 est un quotient d'un groupe de Lie complexe par un r\'eseau cocompact.

\begin{proposition} \label{parallelisable} Soit $M$ une vari\'et\'e complexe  compacte  parall\'elisable  de dimension $n$, quotient d'un groupe de Lie complexe connexe $G$ par un r\'eseau cocompact $\Gamma$ de $G$. Alors~:

(i) $M$ admet des connexions affines holomorphes. L'image r\'eciproque d'une telle connexion sur $G$ est une connexion affine holomorphe invariante par translation \`a droite.

(ii) $M$ admet des connexions affines holomorphes plates.

(iii) $M$ admet une connexion affine holomorphe   plate sans torsion si et seulement s'il existe un morphisme injectif d'alg\`ebres de Lie  $i :  \mathcal G \to gl(n, \CC) \ltimes  \CC^n$ 
dont l'image intersecte trivialement l'isotropie $gl(n, \CC)$, o\`u $\mathcal G$ est l'alg\`ebre de Lie de $G$ et $gl(n, \CC) \ltimes \CC^n$ est l'alg\`ebre de Lie du  groupe affine de $\CC^n$.
\end{proposition}

\begin{corollaire}  \label{semi-simple} Si $G$ est un groupe de Lie complexe semi-simple, alors les  vari\'et\'es complexes compactes parall\'elisables $M=G /  \Gamma$ n'admettent  aucune connexion affine holomorphe
sans torsion plate.
\end{corollaire}

\begin{corollaire} \label{dim3}  Si $M = G/ \Gamma$ est une vari\'et\'e complexe compacte parall\'elisable  de dimension trois, alors $M$ admet une connexion affine holomorphe plate sans torsion si et seulement si $G$ est r\'esoluble.
\end{corollaire}

\section {Connexions affines holomorphes}

Passons \`a la preuve de la proposition~\ref{parallelisable}. 

\begin{demonstration} 

(i) Soient  $X_{1}, \ldots, X_{n}$, des champs de vecteurs invariants par translation \`a droite sur $G$. Toute connexion $\nabla$ invariante par translation \`a droite sur $G$ est caract\'eris\'ee par des constantes $\Gamma^k_{ij} \in \CC$ telles que $\nabla_{X_{i}}X_{j} =\sum_{k} \Gamma^k_{ij} X_{k}$. Une telle connexion descend sur $M$.

Inversement, si $\nabla$ est une connexion affine holomorphe sur $M$,  les coefficients de Christoffel $\Gamma^k_{ij}$  relatifs \`a la famille de champs de vecteurs holomorphes
sur $M$ dont l'image r\'eciproque sur $G$ est $X_{1}, \ldots, X_{n}$, sont des fonctions holomorphes et donc constantes sur la vari\'et\'e compacte $M$.

(ii) Si les constantes $\Gamma^k_{ij}$ sont toutes nulles, la connexion $\nabla$ est plate. Elle se rel\`eve en une connexion plate bi-invariante sur $G$. La torsion $T(X_{i},X_{j})=
\nabla_{X_{i}}X_{j}-\nabla_{X_{j}}X_{i} - \lbrack X_{i}, X_{j} \rbrack$ est nulle si et seulement si $G$ est ab\'elien.

(iii)  Supposons d'abord que $\nabla$ est une connexion holomorphe plate sans torsion sur $G/ \Gamma$. Comme $\nabla$ est localement isomorphe \`a la connexion standard
de $\CC^n$, l'alg\`ebre de Lie du pseudo-groupe des isomorphismes  locaux  de $\nabla$ est  $gl(n, \CC) \ltimes \CC^n$. Par ailleurs, d'apr\`es le point (i), il existe dans l'alg\`ebre de Lie du pseudo-groupe des isomorphismes  locaux  de $\nabla$ une copie de $\mathcal G$ agissant transitivement (et donc intersectant trivialement l'isotropie $gl(n,\CC)$).

R\'eciproquement, si le morphisme $i$ existe, il engendre une action (\`a droite)  affine localement libre du rev\^etement universel  $\tilde G$ de $G$ sur $\CC^n$ telle que l'orbite de
$0$ est ouverte. Ceci munit le groupe de Lie $\tilde G$ d'une structure affine  complexe  invariante par translation \`a droite. Or, $M$ est biholomorphe \`a un quotient $\tilde G / \tilde \Gamma$, o\`u $\tilde \Gamma$
est un r\'eseau cocompact de $\tilde G$. Par cons\'equent, $M$ h\'erite de la structure affine complexe  $\tilde G$-invariante \`a droite de $\tilde G$.
\end{demonstration} 

En d\'eduisons \`a pr\'esent le corollaire~\ref{semi-simple}.

\begin{demonstration}

On appliquera  le point (iii) de la proposition~\ref{parallelisable} au groupe $G$ suppos\'e semi-simple complexe et de dimension $n$.

Consid\'erons un morphisme injectif $i$ de l'alg\`ebre de Lie $\mathcal G$ de $G$ dans $gl(n, \CC) \ltimes  \CC^n$. Comme $\mathcal G$ est semi-simple, la projection $p_1 \circ i$ sur le premier facteur est \'egalement un morphisme injectif (sinon $\mathcal G$ admettrait un id\'eal ab\'elien non trivial). Le lemme de Whitehead~\cite{HS} affirme alors que le premier groupe de cohomologie de $(p_1 \circ i)(\mathcal G)$ \`a coefficients dans la repr\'esentation induite  par la repr\'esentation canonique de $gl(n,\CC)$ sur $\CC^n$ s'annule. Ceci implique que, \`a  automorphisme  interne pr\`es, les seuls morphismes injectifs de $\mathcal G$ dans $gl(n, \CC) \ltimes \CC^n$ sont \`a image dans $gl(n, \CC)$. 

Il reste \`a montrer qu'une repr\'esentation lin\'eaire de $\mathcal G$ sur $\CC^n$ n'admet aucune orbite ouverte. Consid\'erons une telle repr\'esentation et soient $K_{1}, K_{2}, \ldots, K_{n}$ des champs de vecteurs holomorphes sur $\CC^n$ qui sont les champs fondamentaux de l'action de $\mathcal G$ associ\'es \`a une base de $\mathcal G$.

Comme $\mathcal G$ est semi-simple, son repr\'esentation sur $\CC^n$ est \`a valeurs dans l'alg\`ebre de Lie $sl(n, \CC)$ du groupe sp\'ecial lin\'eaire et pr\'eserve la  forme volume
holomorphe $vol=dz_{1} \wedge \ldots \wedge dz_{n}$ sur $\CC^n$.

Supposons par l'absurde qu'une telle repr\'esentation aurait une orbite
ouverte non triviale $\mathcal O$.  Comme $\mathcal G$ est unimodulaire, la fonction holomorphe $vol(K_{1}, K_{2}, \ldots, K_{n})$ est constante (non nulle)  sur $\mathcal O$, et donc 
sur $\CC^n$. Ceci contredit le fait que l'action fixe $0$ et donc tous les $K_{i}$  s'annulent \`a l'origine.
\end{demonstration}

En d\'eduisons le  corollaire~\ref{dim3}.

\begin{demonstration} 

On appliquera  de nouveau  le point (iii) de la proposition~\ref{parallelisable}.

Les alg\`ebres de Lie complexes unimodulaires de dimension trois sont~: $sl(2, \CC)$, $\CC^3$ et les alg\`ebres de Lie (r\'esolubles) $heis$ et $sol$~\cite{Kir}. Rappelons que les alg\`ebres
$heis$ et $sol$ sont engendr\'ees par des g\'en\'erateurs $e_{1}, e_{2}, e_{3}$, avec les relations de crochet  respectivement $\lbrack e_{1}, e_{2} \rbrack =e_{3}, \lbrack e_{1}, e_{3} \rbrack =\lbrack e_{2}, e_{3} \rbrack =0$
et $\lbrack e_{1}, e_{2} \rbrack =e_{2}, \lbrack e_{1}, e_{3} \rbrack =-e_{3}, \lbrack e_{2},
e_{3} \rbrack =0$.

Le cas de l'alg\`ebre de Lie semi-simple $sl(2,\CC)$ vient d'\^etre trait\'e.

Bien s\^ur $\CC^3$ agit librement par translation  sur  $\CC^3$. 

Il reste \`a exhiber dans l'alg\`ebre de Lie du groupe affine de $\CC^3$ des copies de $heis$ ou $sol$ qui intersectent  trivialement l'isotropie. Soit $(f_{1},f_{2},f_{3})$ la base canonique de $\CC^3$.
Les \'el\'ements  $e_{1}=(A,f_{1}), e_{2}=(0, f_{2}),  e_{3}=(0,f_{3}) \in gl(3, \CC) \ltimes \CC^3$ engendrent  une alg\`ebre de Lie  isomorphe \`a $heis$ ou \`a $sol$
selon que  $A \in gl(3,\CC)$ v\'erifie $Af_{2}=f_{3}, Af_{3}=0$,  ou bien $Af_{2}=f_{2}, Af_{3}=-f_{3}$.
\end{demonstration}

Rappelons qu'une connexion affine holomorphe  sans torsion $\nabla$ sur une vari\'et\'e complexe $M$ de dimension $n$ est dite {\it projectivement plate} s'il existe  un atlas de $M$ \`a valeurs dans des ouverts de $P^n(\CC)$ dont chaque carte redresse les g\'eod\'esiques de $\nabla$ sur des droites de $P^n(\CC)$ (sans n\'ecessairement pr\'eserver le  param\`etrage). Dans ce cas les changements de carte sont des transformations projectives et cet atlas munit $M$ d'une {\it structure projective complexe.}~\cite{Gunning}

\begin{proposition} Toute vari\'et\'e complexe compacte parall\'elisable  de dimension trois $M= G / \Gamma$ admet des  connexions affines holomorphes sans torsion projectivement plates.
\end{proposition}

\begin{demonstration} D'apr\`es le corollaire~\ref{dim3}, il reste \`a prouver le r\'esultat pour $G=SL(2,\CC)$. On commence par construire une structure projective complexe
invariante \`a droite sur $SL(2,\CC)$. Remarquons que $P^3(\CC)$ est l'espace projectif sur l'espace vectoriel des polyn\^omes homog\`enes complexes  en deux variables de degr\'e $3$. L'action lin\'eaire (par changement lin\'eaire de variable) de $SL(2,\CC)$ sur cet espace vectoriel se projectivise et admet dans $P^3(\CC)$ l'orbite ouverte qui provient des polyn\^omes homog\`enes qui sont produits de trois formes lin\'eaires distinctes (l'action de $SL(2,\CC)$ sur $P^1(\CC)$ est trois fois transitive). Ceci fournit une action projective
 de $SL(2,\CC)$ sur $P^3(\CC)$ qui poss\`ede une orbite ouverte.  On obtient donc  une structure projective complexe invariante par translation sur $SL(2, \CC)$. Celle-ci descend bien sur $SL(2, \CC) / \Gamma$.

Comme le fibr\'e canonique de $M$ est trivial, cette structure projective complexe correspond \`a une connexion affine holomorphe sans torsion projectivement plate sur $M$ (voir
~\cite{KO}, formule (3.6), pages 78-79).

Une m\'ethode plus directe est de consid\'erer la connexion (holomorphe sans torsion) standard sur $SL(2,\CC)$,  d\'etermin\'ee par $\nabla_xy=\frac{1}{2} \lbrack x,y \rbrack$, pour  $x,y \in sl(2, \CC)$. On v\'erifie que celle-ci est projectivement plate (formellement
les calculs sont les m\^emes que pour la sph\`ere $S^3$).
\end{demonstration}

{\bf Quotients exotiques  de Ghys}.
 En dimension $3$, des exemples  in\'edits de  vari\'et\'es complexes compactes  \'equip\'ees de connexions affines holomorphes  et 
     dont le fibr\'e tangent n'est pas holomorphiquement trivial ont \'et\'e construits dans~\cite{Gh}. Ces exemples 
     s'obtiennent \`a
     partir  de $SL(2,\CC)/ \Gamma$ par d\'eformation  de la structure complexe selon le proc\'ed\'e suivant.
     
     D'apr\`es~\cite{Gh} il existe
 des morphismes de groupe $u : \Gamma \to SL(2,\CC)$ tels que l'action \`a droite de $\Gamma$ sur $SL(2,\CC)$ donn\'ee
 par :
 $$(m,\gamma) \in SL(2,\CC) \times \Gamma \to u(\gamma^{-1})m\gamma \in SL(2,\CC)$$ 
 est libre et totalement discontinue. Le quotient est une vari\'et\'e complexe compacte $M(u,\Gamma)$ qui, en g\'en\'eral,  n'est pas parall\'elisable.
 
 \begin{proposition} Les vari\'et\'es $M(u, \Gamma)$ poss\`edent des connexions affines holomorphes plates et des connexions affines holomorphes  sans torsion projectivement plates, mais aucune connexion affine holomorphe sans torsion plate .
 \end{proposition}
 
 \begin{demonstration} 
 
 D'apr\`es la preuve de la proposition \ref{parallelisable}, $SL(2, \CC)$ admet une connexion affine holomorphe plate bi-invariante. Celle-ci descend en une connexion plate
 sur  $M(u, \Gamma)$.
 
 La connexion standard de $SL(2,\CC)$ est projectivement plate sans torsion et bi-invariante~: elle descend sur $M(u, \Gamma)$.
 
 Consid\'erons maintenant une connexion affine holomorphe quelconque $\nabla$ sur $M(u, \Gamma)$. La diff\'erence entre $\nabla$ et la connexion standard est un tenseur holomorphe~\cite{Gunning}.
 Il est prouv\'e dans~\cite{Gh} que tout
 tenseur holomorphe sur $M(u,\Gamma)$  se rel\`eve en un tenseur holomorphe sur $SL(2,\CC)$ invariant par translation  \`a droite.
 En particulier, $\nabla$ se rel\`eve en une connexion invariante  par translation \`a droite sur $SL(2,\CC)$ et donc l'alg\`ebre de Lie du pseudo-groupe des isomorphismes locaux 
 de $\nabla$ contient une copie de $sl(2,\CC)$ agissant transitivement. On conclut comme dans la preuve du corollaire~\ref{semi-simple} que $\nabla$ n'est pas plate sans torsion.
 \end{demonstration}

\section{Conclusions} 

La question naturelle qui reste ouverte est de savoir si toute vari\'et\'e complexe compacte $M$ qui admet une connexion affine holomorphe poss\`ede \'egalement une connexion
affine holomorphe plate (quitte  \`a supposer \'eventuellement que les classes de Chern $c_{i}(M, \QQ)$ s'annulent).

Dans  le contexte  plus g\'en\'eral des fibr\'es vectoriels holomorphes,  des  r\'esultats dans ce sens sont connus  sur certaines  vari\'et\'es projectives~\cite{Biswas, Simpson} et sur les vari\'et\'es complexes parall\'elisables~\cite{Winkelmann}.

   {\small

{\footnotesize

\vspace{2cm}

Sorin Dumitrescu,
D\'epartement de Math\'ematiques d'Orsay, Bat. 425,
U.M.R.   8628  C.N.R.S., Univ. Paris-Sud (11),
91405 Orsay Cedex,
France

Sorin.Dumitrescu@math.u-psud.fr

 \end{document}